\documentclass[12pt,reqno]{amsart}

\usepackage{amsfonts}
\usepackage{amsmath}
\usepackage{amssymb}

\usepackage{mathrsfs}
\usepackage{hyperref}

\usepackage{epsfig}
\usepackage{graphicx}

\numberwithin{equation}{section}
\theoremstyle{plain}
\newtheorem{thm}{Theorem}[section]

\theoremstyle{definition}

\begin{document}

\title[Initial problem for time-fractional and space-singular equation]
{Non-local initial problem for second order time-fractional and space-singular equation}

\author[Erkinjon Karimov]{Erkinjon Karimov}
\address{
  Erkinjon Karimov:
  \endgraf
  Institute of Mathematics
  \endgraf
  National University of Uzbekistan
  \endgraf
  Tashkent, 100125
  \endgraf
  Uzbekistan
  \endgraf
  {\it E-mail address} {\rm erkinjon@gmail.com}
  }

  \author[Murat Mamchuev]{Murat Mamchuev}
\address{
  Murat Mamchuev:
  \endgraf
  Department of Theoretical and Mathematical Physics
  \endgraf
  Institute of Applied Mathematics and Automation
  \endgraf
  Shortanova str. 89-A, Nalchik, 360000
  \endgraf
  Kabardino-Balkar Republic
  \endgraf
  Russia
  \endgraf
  {\it E-mail address} {\rm mamchuev@rambler.ru}
 }

\author[Michael Ruzhansky]{Michael Ruzhansky}
\address{
  Michael Ruzhansky:
  \endgraf
  Department of Mathematics
  \endgraf
  Imperial College London
  \endgraf
  180 Queen's Gate, London, SW7 2AZ
  \endgraf
  United Kingdom
  \endgraf
  {\it E-mail address} {\rm m.ruzhansky@imperial.ac.uk}
  }

\thanks{The last author was supported in parts by the EPSRC
grant EP/K039407/1 and by the Leverhulme Grant RPG-2014-02.
No new data was collected or generated during the course of research. }

\date{\today}


\begin{abstract}
In this work, we consider an initial problem for second order partial differential equations with Caputo fractional derivatives in the time-variable and Bessel operator in the space-variable. For non-local boundary conditions, we present a solution of this problem in an explicit form representing it by the Fourier-Bessel series. The obtained solution is written in terms of multinomial Mittag-Leffler functions and first kind Bessel functions.
\end{abstract}

\maketitle

\section{Introduction and formulation of a problem}

It is well-known that partial differential equations are playing a key role in constructing mathematical models for many real-life processes. Especially, during the last decades, many applications of various kinds of fractional differential equations became target of intensive research due to both theoretical and practical reasons, see e.g. \cite{KST06} for an overview. Many kinds of boundary problems, including direct \cite{nakh} and inverse problems \cite{isak}, were formulated for different type of PDEs of integer order including several differential operators of fractional order.

We note works \cite{gejbhal}-\cite{gorluch} devoted to studying partial differential equations with multiple Caputo derivatives. Precisely, in \cite{gorluch} the authors studied fractional differential equations with Caputo fractional derivatives and using the operational method, solutions of initial boundary problem for those equations were obtained in an explicit form involving a multinomial Mittag-Leffler function. Certain properties of this function were obtained by Li, Liu and Yamamoto \cite{yamli} and applied to studying initial-boundary problems for time-fractional diffusion equations with positive constant coefficients. Later, Liu \cite{liu} established strong maximum principle for fractional diffusion equations with multiple Caputo derivatives and investigated a related inverse problem. Daftardar-Gejji and Bhalikar \cite{gejbhal}, using  the method of separation of variables, solved some boundary-value problems for multi-term fractional diffusion-wave equation.

We also note works related to the Bessel operator. In \cite{masmes}, the initial inverse problem for the heat equation with Bessel operator was investigated. Inverse initial and inverse source problems for time-fractional diffusion equation with zero order Bessel operator were recently studied in \cite{fatma}. Direct and inverse problems for PDEs containing two-term time fractional Caputo derivatives of orders up to 1, and Bessel operator of order $\nu$ were investigated in \cite{akmr}.

The consideration of non-local initial conditions is often justified by practical usage in certain real-life processes. For instance, when the initial temperature for the heat equation is not given instantly, but there is an information related with the temperature on a certain time interval that can be described by a non-local initial condition in a simple form. Boundary-value problems with non-local initial conditions were considered in works \cite{pao} for reaction-diffusion equations, in \cite{shop} for heat equation, in \cite{rk1}-\cite{rk2} for degenerate parabolic equations, and in \cite{ker} for a mixed parabolic equation.

For integer orders much more is known, and for a review of different questions of time decay of solutions for hyperbolic equations with integer order derivatives we can refer to \cite{RS}.

In the present work we deal with the non-local initial boundary problem for multi-term time fractional PDE with Bessel operator of order $\nu$. We use Fourier-Bessel series expansion in order to find the explicit solution for the considered problem, yielding also its existence. Because of the singularities in the Bessel operator such conditions appear naturally in space variables.

We note that most of the current literature deals with diffusion type equations considering time-derivatives of orders up to 1, see e.g. \cite{yamli}. One of the novelties of the present paper is that we consider wave type equations allowing fractional derivatives of order up to 2, with additional fractional dissipation type terms. If there are multiple fractional time-derivative terms, a multinomial Mittag-Leffler function appears in the representation of solutions.

\medskip
Let us now describe the problem in more detail. We consider the equation
\begin{equation}\label{eq1}
L(u)-B_\nu(u)=f(t,x)
\end{equation}
in a rectangular domain $D=\left\{(x,t):\, 0<x<1,\,\,0<t<T\right\}$, $T>0$,
where $f(t,x)$ is a given function,
\begin{equation}\label{eq2-3}
L(u)=\partial_{0t}^{\alpha}u(t,x)-\sum\limits_{i=1}^n \lambda_i \partial_{0t}^{\alpha_i}u(t,x)
\end{equation}
is the time component of the equation, with orders 
$$0<\alpha_i\leq1,\quad \alpha_i\leq \alpha\le 2,\quad n\in \mathbb{N}, \quad \lambda_i\in \mathbb{R},$$ 
and
\begin{equation}\label{eq2-3b}
  B_\nu(u)=u_{xx}(t,x)+\frac{1}{x}u_x(t,x)-\frac{\nu^2}{x^2}u(t,x)
\end{equation}
is the Bessel part of the equation with $\nu>0$. Here
\begin{equation*}
\partial_{0t}^{\beta}g(t)=\left\{ \begin{aligned}
  & \frac{1}{\Gamma \left( k-\alpha  \right)}\int\limits_{0}^{t}{\frac{{g}^{(k)}\left( z \right)}{{{\left( t-z \right)}^{\alpha-k+1 }}}dz,\,}\,\,\alpha\notin \mathbb{N}_0, \\
 & \frac{d^kg(t)}{dt^k},\quad\,\,\,\,\,\,\,\,\,\,\,\,\,\,\,\,\,\,\,\,\,\,\,\,\,\,\,\,\,\,\,\,\,\,\,\,\,\,\,\,\,\,\,\,\alpha \in\mathbb{N}, \\
& g(t), \quad\quad\,\,\,\,\,\,\,\,\,\,\,\,\,\,\,\,\,\,\,\,\,\,\,\,\,\,\,\,\,\,\,\,\,\,\,\,\,\,\,\,\,\,\,\,\,\alpha=0,
\end{aligned} \right.
\end{equation*}
is a fractional differential operator of Caputo type, where $k=[\alpha]+1$, and $[\alpha]$ is the integer part of $\alpha$. We can refer to \cite{KST06} for further details on the Caputo fractional derivative operators.

\medskip
The non-local initial boundary problem for equation \eqref{eq1}-\eqref{eq2-3b} is formulated as follows:

\bigskip
{\bf Problem: % $NI_\alpha$.
} Let $M\in\mathbb R$. To find a solution $u(t,x)$ of equation \eqref{eq1}-\eqref{eq2-3b} in $D$, which satisfies
\begin{itemize}
  \item[(i)] regularity conditions $u\in W$ with
  \begin{equation}\label{eq4}
W=\left\{u(t,x):\, u\in C(\overline{D}), \,\,u_{xx},\,\, \partial_{0t}^{\alpha}u\in C(D),\,\,\int\limits_0^1 \sqrt{x} |u(t,x)|dx<+\infty\right\};
  \end{equation}
  \item[(ii)] boundary and non-local initial conditions
  \begin{eqnarray}\label{eq5-6}
      \lim\limits_{x\rightarrow 0} x u_x(t,x)=0,\,\, u(t,1)=0, \\ \label{eq5-6b}
    u(0,x)+M u(T,x)=0,\,\, 0\le x\le 1,\,\,\,\, [\alpha]\cdot u_t(0,t)=0,\,\,0<x<1.
      \end{eqnarray}
\end{itemize}

\section{Main result}

The main result of this note is the following well-posedness theorem for the initial problem
\eqref{eq1}-\eqref{eq5-6b}. The interesting part are the conditions on $f$ allowing one to handle the singularities in the coefficients of the Bessel operator, and the non-resonance conditions \eqref{EQ:M} relating the parameter $M$ with the length $T$ of the time interval, coefficients and fractional orders of time-derivatives, through the multinomial Mittag-Leffler function.

\begin{thm}\label{THM:main}
Assume that
\begin{itemize}
\item $f(x,t)$ is differentiable four times with respect to $x$;
\item $f(0,t)=f'(0,t)=f''(0,t)=f'''(0,t)=0,\,\,f(1,t)=f'(1,t)=f''(1,t)=0$;
\item $\frac{\partial^4 f(x,t)}{\partial x}$ is bounded;
\item $f(x,t)$ is continuous and continuously differentiable with respect to $t$,
\end{itemize}
and non-resonance conditions
\begin{equation}\label{EQ:M}
M\neq -\frac{1}{E_{(\alpha-\alpha_1,\alpha-\alpha_2,...,\alpha),1}(\lambda_1 T^{\alpha-\alpha_1},...,\lambda_n T^{\alpha-\alpha_n},-\gamma_k^2T^\alpha)}
\end{equation}
hold for all $k=1,2,\ldots$. Then there exists a unique solution of the problem \eqref{eq1}-\eqref{eq5-6b}, and it can be written in the following form:

\begin{equation}\label{EQ:sol}
\begin{aligned}
&u(t,x)=\sum\limits_{k=1}^\infty \left[\int\limits_0^t z^{\alpha-1}E_{(\alpha-\alpha_1,\alpha-\alpha_2,...,\alpha-\alpha_n,\alpha),\alpha}(\lambda_1 z^{\alpha-\alpha_1},...,\lambda_n z^{\alpha-\alpha_n},-\gamma_k^2 z^\alpha)f_k(t-z)dz-\right.\\
&-\frac{M}{1+ME_{(\alpha-\alpha_1,\alpha-\alpha_2,...,\alpha),1}(\lambda_1 T^{\alpha-\alpha_1},...,\lambda_n T^{\alpha-\alpha_n},-\gamma_k^2T^\alpha)}\times\\
&\times \int\limits_0^T z^{\alpha-1}E_{(\alpha-\alpha_1,\alpha-\alpha_2,...,\alpha-\alpha_n,\alpha),\alpha}(\lambda_1 z^{\alpha-\alpha_1},...,\lambda_n z^{\alpha-\alpha_n},-\gamma_k^2 z^\alpha)f_k(T-z)dz+\\
&\left.+E_{(\alpha-\alpha_1,\alpha-\alpha_2,...,\alpha),1}(\lambda_1 t^{\alpha-\alpha_1},...,\lambda_n t^{\alpha-\alpha_n},-\gamma_k^2t^\alpha)\right]J_\nu(\gamma_k x).
\end{aligned}
\end{equation}
\end{thm}

The numbers $\gamma_k$ and the functions appearing in the formula \eqref{EQ:sol} are explained in Section \ref{SEC:series}. Briefly, here $E_{(\alpha-\alpha_1,\alpha-\alpha_2,...,\alpha-\alpha_n,\alpha),\alpha}(\cdot)$ is the multinomial Mittag-Leffler function, $J_\nu$ are the first kind Bessel functions,
$\gamma_k$ are their positive zeros, $\lambda_k$ are coefficients in the operator \eqref{eq2-3}, and $f_k$ are Bessel expansions of $f$.

For the proof of Theorem \ref{THM:main} we start finding a formal solution in a series form and the convergence of the appearing series will be shown in Section \ref{SEC:proof}.

\subsection{Representation of a solution}
\label{SEC:series}

We start by a formal discussion of the representation of solutions in \eqref{EQ:sol}.
Let
\begin{equation*}
J_\nu(z)=\sum\limits_{i=0}^\infty \frac{(-1)^i(z/2)^{2i+\nu}}{i!(i+\nu)!}
\end{equation*}
be the Bessel function of the first kind (see e.g. \cite{Wats}).  It is known that for $\nu>0$, the Bessel function $J_\nu(z)$ has countably many zeros, moreover, they are real and have pairwise opposite signs. Denote the $k^{th}$ positive root of the equation $J_\nu(z)=0$ by $\gamma_k$, $k=1,2,\ldots$.
%, for which nontrivial solutions of the eigenvalue problem exist.
For large $k$, we have (see \cite{tols})
\begin{equation*}
\gamma_k\simeq k\pi+\frac{\nu\pi}{2}-\frac{\pi}{4}.
\end{equation*}

We now expand functions $u(t,x)$ and $f(t,x)$ in the Fourier-Bessel series (see e.g. \cite{tols}), writing them in the form

\begin{equation}\label{eq7}
u(t,x)=\sum\limits_{k=1}^\infty U_k(t) J_\nu\left(\gamma_k x\right),
\end{equation}
\begin{equation}\label{eq8}
f(t,x)=\sum\limits_{k=1}^\infty f_k(t) J_\nu\left(\gamma_k x\right),
\end{equation}
where
\begin{equation}\label{eq9}
U_k(t)=\frac{2}{J_{\nu+1}^2(\gamma_k)}\int\limits_0^1 u(x,t)\, xJ_\nu\left(\gamma_k x\right)dx,
\end{equation}
\begin{equation}\label{eq10}
f_k(t)=\frac{2}{J_{\nu+1}^2(\gamma_k)}\int\limits_0^1 f(x,t)\, xJ_\nu\left(\gamma_k x\right)dx.
\end{equation}

It is known that if a function $g=g(x)$ is piecewise continuous on  $[0,l]$ and satisfies
$$\int\limits_0^l\sqrt{x}|g(x)|dx<+\infty,$$
then for $\nu>-1/2$, the Fourier-Bessel series converges at every point $x_0\in (0,l)$, see e.g. \cite{tols}. Since we are looking for a function $u(t,x)\in W$, it satisfies these required conditions in order to be represented by a Fourier-Bessel series.

We substitute \eqref{eq7}-\eqref{eq8} into equation \eqref{eq1} and obtain the eigenvalue problem
\begin{equation}\label{eq11}
L(U_k)+\gamma_k^2U_k(t)=f_k(t).
\end{equation}
According to \cite{gorluch}, the solution for equation \eqref{eq11} satisfying initial conditions
\begin{equation}\label{eq12}
U_k(0)=A,\,\,[\alpha] U_k'(0)=0,
\end{equation}
can be represented in the form
\begin{multline}\label{eq13}
U_k(t)= 
\int\limits_0^t z^{\alpha-1}E_{(\alpha-\alpha_1,\alpha-\alpha_2,...,\alpha-\alpha_n,\alpha),\alpha}(\lambda_1 z^{\alpha-\alpha_1},...,\lambda_n z^{\alpha-\alpha_n},-\gamma_k^2 z^\alpha)\times \\ \times f_k(t-z)dz+A \overline{U_0}(t),
\end{multline}
where
\begin{multline}\label{eq14}
\overline{U_0}(t)=1+\sum\limits_{i=1}^n \lambda_it^{\alpha-\alpha_i}E_{(\alpha-\alpha_1,\alpha-\alpha_2,...,\alpha-\alpha_n,\alpha),1+\alpha-\alpha_i}(\lambda_1 t^{\alpha-\alpha_1},...,\lambda_n t^{\alpha-\alpha_n},-\gamma_k^2 t^\alpha)-\\
-\gamma_k^2t^\alpha E_{(\alpha-\alpha_1,\alpha-\alpha_2,...,\alpha-\alpha_n,\alpha),1+\alpha}(\lambda_1 t^{\alpha-\alpha_1},...,\lambda_n t^{\alpha-\alpha_n},-\gamma_k^2 t^\alpha),\end{multline}
and
\begin{equation}\label{eq15}
E_{(a_1,a_2,...,a_n),b}(z_1,z_2,...,z_n)=\sum\limits_{k=0}^\infty \sum\limits_{\begin{array}{l}l_1+l_2+...+l_n=k\\
l_1\geq0,...,l_n\geq 0\\
\end{array}}\frac{k!}{l_1!...l_n!}\frac{\prod\limits_{i=1}^n z_i^{l_i}}{\Gamma(b+\sum\limits_{i=1}^n a_i l_i)}
\end{equation}
is the multinomial Mittag-Leffler function (\cite{gorluch}).
From \eqref{eq13} we find that
\begin{multline}\label{eq16}
U_k(T)= 
\int\limits_0^T z^{\alpha-1}E_{(\alpha-\alpha_1,\alpha-\alpha_2,...,\alpha-\alpha_n,\alpha),\alpha}(\lambda_1 z^{\alpha-\alpha_1},...,\lambda_n z^{\alpha-\alpha_n},-\gamma_k^2 z^\alpha)\times \\ \times f_k(T-z)dz+A \overline{U_0}(T).
\end{multline}
Considering $U_k(0)=A$, from the first condition in \eqref{eq5-6}, we get the relation
\begin{equation*}
A+MU_k(T)=0.
\end{equation*}
Using \eqref{eq16} we find $A$ to be
\begin{multline*}
A=-\frac{M}{1+M\overline{U_0}(T)}\int\limits_0^T z^{\alpha-1}E_{(\alpha-\alpha_1,\alpha-\alpha_2,...,\alpha-\alpha_n,\alpha),\alpha}(\lambda_1 z^{\alpha-\alpha_1},...,\lambda_n z^{\alpha-\alpha_n},-\gamma_k^2 z^\alpha)\times\\
\times f_k(T-z)dz.
\end{multline*}
Substituting this value of $A$ into \eqref{eq13}, we rewrite it as
\begin{multline}\label{eq17}
U_k(t)=\int\limits_0^t z^{\alpha-1}E_{(\alpha-\alpha_1,\alpha-\alpha_2,...,\alpha-\alpha_n,\alpha),\alpha}(\lambda_1 z^{\alpha-\alpha_1},...,\lambda_n z^{\alpha-\alpha_n},-\gamma_k^2 z^\alpha)f_k(t-z)dz-\\
-\frac{M\overline{U_0}(t)}{1+M\overline{U_0}(T)}\int\limits_0^T z^{\alpha-1}E_{(\alpha-\alpha_1,\alpha-\alpha_2,...,\alpha-\alpha_n,\alpha),\alpha}(\lambda_1 z^{\alpha-\alpha_1},...,\lambda_n z^{\alpha-\alpha_n},-\gamma_k^2 z^\alpha)\times \\
\times f_k(T-z)dz.
\end{multline}
If we use the formula (see \cite{liu})
\begin{equation*}
\begin{aligned}
&1+\sum\limits_{j=1}^{n+1}\lambda_jt^{\alpha-\alpha_j}E_{(\alpha-\alpha_1,\alpha-\alpha_2,...,\alpha-\alpha_{n+1}),1+\alpha-\alpha_j}(\lambda_1 t^{\alpha-\alpha_1},...,\lambda_{n+1} t^{\alpha-\alpha_{n+1}})\\
&=E_{(\alpha-\alpha_1,\alpha-\alpha_2,...,\alpha-\alpha_{n+1}),1}(\lambda_1 t^{\alpha-\alpha_1},...,\lambda_{n+1} t^{\alpha-\alpha_{n+1}}),
\end{aligned}
\end{equation*}
representation \eqref{eq14} can be rewritten as
\begin{equation*}
\overline{U_0}(t)=E_{(\alpha-\alpha_1,\alpha-\alpha_2,...,\alpha-\alpha_{n+1}),1}(\lambda_1 t^{\alpha-\alpha_1},...,\lambda_n t^{\alpha-\alpha_n},-\gamma_k^2 t^\alpha).
\end{equation*}
Denoting
\begin{multline}\label{eq18-19}
F_k(t)=
\int\limits_0^t z^{\alpha-1}E_{(\alpha-\alpha_1,\alpha-\alpha_2,...,\alpha-\alpha_n,\alpha),\alpha}(\lambda_1 z^{\alpha-\alpha_1},...,\lambda_n z^{\alpha-\alpha_n},-\gamma_k^2 z^\alpha)f_k(t-z)dz,
\\
\overline{U_0}(t)=E_{(\alpha-\alpha_1,\alpha-\alpha_2,...,\alpha),1}(\lambda_1 t^{\alpha-\alpha_1},...,\lambda_n t^{\alpha-\alpha_n},-\gamma_k^2t^\alpha),
\end{multline}
we rewrite the function \eqref{eq17} as
\begin{equation}\label{eq20}
U_k(t)=F_k(t)-\frac{M}{1+M\overline{U_0}(T)}F_k(T)\overline{U_0}(t).
\end{equation}
We note that the above expression is well-defined in view of the non-resonance conditions \eqref{EQ:M}, that is
\begin{equation*}
M\neq -\frac{1}{\overline{U_0}(T)}=-\frac{1}{E_{(\alpha-\alpha_1,\alpha-\alpha_2,...,\alpha),1}(\lambda_1 T^{\alpha-\alpha_1},...,\lambda_n T^{\alpha-\alpha_n},-\gamma_k^2T^\alpha)}
\end{equation*}
holds for all $k$.
Finally, based on \eqref{eq20}, we rewrite our formal solution as

\begin{equation}\label{eq21}
u(t,x)=\sum\limits_{k=1}^\infty \left[F_k(t)-\frac{M}{1+M\overline{U_0}(T)}F_k(T)\overline{U_0}(t)\right]J_\nu(\gamma_kx).
\end{equation}

\subsection{Justification of formal solution}
\label{SEC:proof}

In this section we prove convergence of the obtained infinite series corresponding to functions $u(t,x),\,u_{xx}(t,x)$ and $\partial_{0t}^\alpha u(t,x)$.

In order to prove the convergence of these series, we use the estimate for the Mittag-Leffler function \eqref{eq15}, obtained in \cite[Lemma 3.2]{yamli}, of the form

\begin{equation}\label{eq22}
|E_{(\alpha-\alpha_1,\alpha-\alpha_2,...,\alpha-\alpha_n),\rho}(z_1,z_2,...,z_n)|\leq \frac{C}{1+|z_1|}.
\end{equation}

Let us first prove the convergence of series \eqref{eq21}. For this, we collect several other known estimates. First, we use the  following theorem on the estimate of the Fourier-Bessel coefficient:

\begin{thm}[{\cite[p. 231]{tols}}]\label{THM:1}
Let $f(x)$ be a function defined on the interval $[0,1]$ such that $f(x)$ is differentiable $2s$ times  $(s\in\mathbb N)$ and such that

(1) $f(0)=f'(0)=...=f^{(2s-1)}(0)=0$;

(2) $f^{(2s)}(x)$ is bounded (this derivative may not exist at certain points);

(3) $f(1)=f'(1)=...=f^{(2s-2)}(1)=0$.

Then the following inequality is  satisfied by the Fourier-Bessel coefficients of $f(x)$:
\begin{equation}\label{eq23}
|f_k|\leq \frac{c}{\gamma_k^{2s-1/2}}\,\,\,\,(c=const),
\end{equation}
where $\gamma_k$ is the $k^{th}$ positive zero of of the function $J_\nu(x)$.

\end{thm}

In our case we have
$$
f_k(t)=F_k(t)-\frac{M}{1+M\overline{U_0}(T)}F_k(T)\overline{U_0}(t).
$$
According to \eqref{eq18-19}, using estimate \eqref{eq22} and imposing conditions (1)-(3) of Theorem \eqref{THM:1} in the case $s=1$ on the function $f(t,x)$, we get
\begin{equation}\label{eq24}
\left|F_k(t)-\frac{M}{1+M\overline{U_0}(T)}F_k(T)\overline{U_0}(t)\right|\leq \frac{C_1}{\gamma_k^{3/2}}.
\end{equation}

The series \eqref{eq20} then converges absolutely and uniformly on $[0,1]$ in view of the following theorem:

\begin{thm}[{\cite[p. 225]{tols}}]\label{THM:2}
If $\nu\geq 0$, $C>0$, and if the constants $c_k$ satisfy
$$
|c_k|\leq \frac{C}{\gamma_k^{1+\varepsilon}},
$$
for some $\varepsilon>0$, then the series
$$
c_1J_\nu(\gamma_1x)+c_2J_\nu(\gamma_2x)+c_3J_\nu(\gamma_3x)+...
$$
converges absolutely and uniformly on $[0,1]$.
\end{thm}

According to \cite[Theorem 2, p. 236]{tols}, sufficient condition for differentiating the series \eqref{eq7} twice term by term, i.e. for the validity of the relation
\begin{equation}\label{EQ:diff}
u_{xx}(t,x)=\sum\limits_{k=1}^\infty U_k(t)\gamma_k^2 J_\nu''(\gamma_kx)
\end{equation}
will be 
\begin{equation}\label{EQ:quest}
|U_k(t)|\leq \frac{C_2}{\gamma_k^{5/2+\epsilon}},
\end{equation}
where $U_k$ are as in \eqref{eq20}, and $C_2$ and $\epsilon$ are positive constants.
Based on this estimation and considering also (see \cite[p. 233]{tols})
$$
|J_\nu(\gamma_kx)|\leq \frac{C_3}{\sqrt{\gamma_k x}},\quad (C_3=const),
$$
we deduce (see e.g. \cite[p. 236]{tols})
$$
\left|U_k(t)\gamma_k^2J_\nu''(\gamma_kx)\right|\leq \frac{C_2}{\gamma_k^{2+\epsilon}x^2}+\frac{C_3(|\nu|+\nu^2)}{\gamma_k^{3+\epsilon}x^2\sqrt{x}}+\frac{C_4}{\gamma_k^{1+\epsilon}\sqrt{x}},
$$
which provides the convergence of the series \eqref{EQ:diff}.
We note that if we impose conditions on the given function $f(x,t)$ of the form
\begin{itemize}
\item $f(x,t)$ is differentiable four times with respect to $x$;
\item $f(0,t)=f'(0,t)=f''(0,t)=f'''(0,t)=0,\,\,f(1,t)=f'(1,t)=f''(1,t)=0$;
\item $\frac{\partial^4 f(x,t)}{\partial x}$ is bounded,
\end{itemize}
one can see that Theorem \ref{THM:2} implies estimation \eqref{EQ:quest}.

The convergence of series corresponding to $\partial _{0t}^\alpha u(x,t)$, $u_x(x,t)$ can be shown in similar ways, completing the proof of Theorem \ref{THM:main}.

\end{document}